\documentclass[psfig,reqno,a4paper,11pt]{amsart}
\usepackage{amsfonts}

\theoremstyle{plain}
\begingroup

\newtheorem{theorem}{Theorem}[section]
\newtheorem{lemma}[theorem]{Lemma}
\newtheorem{proposition}[theorem]{Proposition}

\newtheorem{definition}[theorem]{Definition}

\theoremstyle{definition}
\theoremstyle{remark}
\numberwithin{equation}{section}
\newcommand{\var}{\varphi}
\newcommand{\e}{\varepsilon}
\newcommand{\Om}{\Omega}

\newcommand{\dx}{\,dx}

\newcommand{\intcauchy}{\mskip 3mu -\mskip -19mu \int}

\newcommand{\R}{{\mathbb R}}

\newcommand{\salt}{\noalign{\vskip .2truecm}}
\newcommand{\parent}[3]{\left #1 {#3} \right #2} 
\newcommand{\graffe}[1]{\parent \{ \}{#1}} 
\newcommand{\bord}{{\partial\Om}}

\input psfig.sty
\newcommand{\bfig}[2]{\begin{figure}\begin{center}\begin{picture}(341.8,#2)(
#1,0)}
\newcommand{\efig}[2]{\end{picture}\caption{#2.}\lbl{#1}\end{center}
\end{figure}}
 
\title[Stability of some unilateral free-discontinuity problems]{Stability of some unilateral free-discontinuity\\ 
  problems in two-dimensional domains}
\author[Fran\c{c}ois Ebobisse]{Fran\c{c}ois Ebobisse}
\address[Fran\c{c}ois Ebobisse]{S.I.S.S.A., Via Beirut 2-4, 34014,
Trieste,
Italy}
\email{ebobisse@sissa.it}
\author[Marcello Ponsiglione]{Marcello Ponsiglione}
\address[Marcello Ponsiglione]{S.I.S.S.A., Via Beirut 2-4, 34014,
Trieste, Italy}
\email{ponsigli@sissa.it}

\begin{document}
\baselineskip3.3ex

\vskip .2truecm
\begin{abstract}
\small{The purpose of this paper is to study
the stability of  some unilateral free-discontinuity problems in two-dimensional domains, 
with the density of the volume part having $p$-growth, with $1<p<\infty$,  under perturbations of the discontinuity sets in the 
Hausdorff metric.
}
\end{abstract}
\maketitle
{\small \tableofcontents}
%\vskip.1truecm
\section{Introduction}\label{Intro}
 Let $\Om\subset\R^2$ be a bounded connected open set with Lipschitz continuous 
boundary $ \partial\Om$ and  let   
 $\partial_D\Om\subset\partial\Om$ be a (non-empty) relatively open subset of 
$\bord$ composed of a finite number of connected components. Let $g\in W^{1,p}(\Om)$
 and let  $f:\Om\times\R^2\to\R$  be a Borel function which satisfies the assumptions (\ref{f1})-(\ref{f2}) below.
 We consider pairs $(u,K)$ with $K$ a compact subset of $\overline\Om$ and 
$u\in L^{1,p}(\Om\setminus K):=\{ v\in L^p_{loc}(\Om\setminus K),\, \nabla v\in L^p(\Om\setminus K,\R^2)\} $ with
 $u=g$ on $\partial_D\Om\setminus K$, which satisfy the following unilateral minimality condition:
\begin{equation}\label{minimalim}
\int_{\Om\setminus K}f(x,\nabla u)\dx\,+\,{\mathcal H}^1(K)\,
\leq\,\int_{\Om\setminus H}f(x,\nabla w)\dx\,+\,{\mathcal H}^1(H),
\end{equation}
 among all compact subsets $H$ of $\overline\Om$ with $H\supset K$ 
and all functions $ w\in L^{1,p}(\Om\setminus H)$ with $w=g$ on $\partial_D\Om\setminus H$. 

Our goal in this paper is to study the stability of the problem 
(\ref{minimalim}) under variations of the compact set $K$ in the 
Hausdorff metric and of the boundary datum $g$ in the strong 
 topology of $W^{1,p}(\Om)$. Precisely, let $(K_h)$ be a 
sequence of compact subsets of $\overline\Om$ which converges to 
 a compact set $K$ in the Hausdorff metric and let 
$(g_h)\subset W^{1,p}(\Om) $ be a sequence which converges 
 strongly to a function $g$ in $W^{1,p}(\Om)$. 
Let $u_h$ be such that the pair $(u_h,K_h)$ is a solution of (\ref{minimalim}) relative to the 
boundary data $g_h$. We are studying the conditions under which the sequence $(\nabla u_h)$ 
converges strongly to $\nabla u$ in $L^p(\Om,\R^2)$ for some function $u$ such that the pair $(u,K)$ 
is a solution of (\ref{minimalim}) relative to the boundary data $g$.

Minimization problems of the type (\ref{minimalim}) arise for instance in the 
mathematical formulation of the irreversible quasi-static growth of brittle fractures 
 based on Griffith's theory of crack growth.
 In this model, the crack path is determined by the competition between
 bulk and surface energy.
The variational  model 
%for the irreversible quasi-static growth of brittle fractures based on  Griffith's theory has been 
 proposed by G.A. Francfort and J.-J. Marigo \cite{FM} is described as follows:  from an initial crack $K_0$ (possibly an empty set), 
 the crack $K(t)$ at a given time $t$ corresponding to a loading $g(t)$ applied to
 $\partial_D\Om$, will minimize the total energy (bulk energy + surface energy) among all the possible cracks $K$ which contain
 the previous one $K(s)$, $s<t$. This  continuum evolution of the cracks 
 during the loading process is obtained as a limit of a discretized evolution described as 
a step by step unilateral minimization problem of the type (\ref{minimalim}).

The precise mathematical formulation of this model has been studied by G. Dal Maso and R. Toader \cite{dmro1, dmro2}
in the special case of linearized elasticity for anti-plane shear 
and for an {\it a priori bound} on the number of connected 
 component of the test cracks. In this case the reference configuration 
is an infinite cylinder $\Om\times\R$, with $\Om\subset\R^2$, and the displacement 
field has the form $v:=(0,0,u)$ where $u$ is a scalar function defined on $\Om$. 
The cracks are assumed also to be of the form $K\times\R$, 
where $K$ is a compact subset of $\overline\Om$.

 Recently, a weak formulation for the variational model of fracture 
growth in the framework of $SBV$ space of special functions of bounded variation,
 has been proposed by  G.A. Francfort and C.J. Larsen \cite{FL} 
for anti-plane shear in higher dimensions. This approach is more 
natural since it is performed in any dimension and with 
 no restrictions on the test cracks. However, the strong formulation in  \cite{dmro1}, 
 based on the Hausdorff convergence of compact sets,  is more 
handable and elementary in two dimensions and leads to the convergence 
in the Hausdorff metric of the cracks obtained in the discretized evolution.

One of the key points in \cite{dmro1} is the stability of (\ref{minimalim}) for $f(x,\xi)=|\xi|^2$, 
which follows from  the stability of the following minimization problem:
\begin{equation}\label{eq0}
\displaystyle\min_{v}\graffe{\int_{\Om\setminus K}f(x,\nabla v)\dx\,\mbox{: }v\in 
L^{1,p}(\Om\setminus K)\,,\quad v=g\,\mbox{ on }\,\partial_D\Om\setminus K}.
\end{equation} 
 Actually the stability of (\ref{eq0}) holds for every $p\leq 2$ under 
the hypotheses of \cite{dmro1} (see \cite{dmebpo}), 
 while in the case $p> 2$ some counter-examples have been given in \cite{dmebpo} and in \cite{ebpo}. 
 The strategy to get  the stability of problem (\ref{minimalim}) 
for every $1<p<\infty$ is to obtain the stability of (\ref{eq0}) 
 using the unilateral minimality condition.

The obstruction to the stability of (\ref{eq0}) when
$p> 2$ is due to the fact that two  connected components of the approximating sequence $(K_h)$
 can approach and touch each other in the limit fracture $K$, leading then to the appearance of 
 a transmission term  in the limit problem. To avoid such phenomena 
 we joint these two connected components by curves of infinitesimal length, obtaining then a new 
 sequence of cracks $(H_h)$ having the properties that $K_h\subset H_h$, $H_h$ converges to $K$, 
$\mathcal H^1(H_h\setminus K_h)\to 0$ and 
any connected component of $H_h$ converges to a connected 
component of the limit fracture $K$.  Then the stability of (\ref{eq0}) along this 
new sequence of cracks $(H_h)$ will follow from Proposition \ref{lemdisjoint}. 
 Now, using the unilateral constraint, 
we obtain the stability of (\ref{eq0}) also along the original sequence of cracks $(K_h)$.
 
We prove our main results (see Theorems  \ref{unistabi} and \ref{genecase}) following the duality approach, i.e., 
through the conjugates (see Section 3), performed in 
  \cite{bv2}, \cite{dmro1} for linear problems, and extended recently in \cite{dmebpo} to nonlinear problems.

\section{Notation and preliminaries}\label{Notprel}
Let $\Om $ be a bounded connected open subset of $\R^2$  
 with Lipschitz continuous boundary $ \partial\Om$. Let   
 $\partial_D\Om\subset\partial\Om$ be a (non-empty) relatively open subset of 
$\bord$ composed of a finite number of connected components and 
$\partial_N\Om:=\bord\setminus\partial_D\Om$. 

Let ${\mathcal K}(\overline\Om)$
be the class of compact subsets of $\overline\Om$ 
and  ${\mathcal K}_m(\overline\Om)$ be the subset of ${\mathcal K}(\overline\Om)$ 
whose elements have at most $m$ connected components. 
We denote $\mathcal K^f_m(\overline \Om)$ the subclass of ${\mathcal
K}_m(\overline\Om)$ whose elements have finite  one-dimensional Hausdorff measure $\mathcal H^1$. 
For every $\lambda >0$, ${\mathcal K}_m^\lambda (\overline\Om)$ 
denotes the class of sets $K$ in ${\mathcal
K}_m(\overline\Om)$ such that $\mathcal H^1(K)\leq\lambda$. 

For any $x\in \Om $ and $\rho>0$, $B(x,\rho )$ denotes 
  the open ball of  $\R^2 $ centered at $x$ with radius $\rho $. 
For any subset $E$ of $\R^2$, $1_E$ is the characteristic function of $E$, $E^c$ is the complement of $E$, and 
$|E|$ is the Lebesgue measure of $E$. Throughout the paper 
$p$ and $q$ are real numbers, with $1<p,\,q<+\infty$ and $p^{-1}+q^{-1}=1$.
\subsection{Deny-Lions spaces}
 Given an open subset
$U$ of $\R^2$, the Deny-Lions space is defined by
$$L^{1,p}(U):=\{u\in L^p_{\rm loc}(U):\, \nabla u\in L^p(U,\R^2)\}.$$
It is well-known that $L^{1,p}(U)$ coincides with the Sobolev space $W^{1,p}(U)$ 
whenever $U$ is bounded and has a Lipschitz continuous boundary. 
It is also known that the set $\{\nabla u:\, u\in L^{1,p}(U)\}$ is a closed subspace of
$L^p(U,\R^2)$. The Deny-Lions spaces $L^{1,p}$ are usually involved in minimization
 problems of the type (\ref{eq1}) below,  in non smooth domains
 where Poincar\'e inequalities  do not hold in general. For further properties of the spaces $L^{1,p}$ we refer the
reader to \cite{deli} and \cite{maz1}.
\subsection{The minimization problem} 
Let $f:\Om\times\R^2\to\R$ be a Borel function which satisfies the following assumptions: there exist 
positive constants $\alpha$, $\beta$, $\gamma$ such that, for almost every $x\in\Om$ and for every $\xi\in\R^2$
\begin{eqnarray}
&\alpha|\xi|^p\leq f(x,\xi)\leq\beta |\xi|^p+\gamma;\label{f1}\\
\salt
&f(x,\cdot)\mbox{ is strictly convex and is of class } C^1.\label{f2}
\end{eqnarray}
Given $K\in {\mathcal
K}(\overline\Om)$ and a function $g\in W^{1,p}(\Om)$, we consider the following  minimization problem
\begin{equation}\label{eq1}
\displaystyle\min_{v}\graffe{\int_{\Om\setminus K}f(x,\nabla v)\dx\,\mbox{: }v\in 
L^{1,p}(\Om\setminus K)\,,\quad v=g\,\mbox{ on }\,\partial_D\Om\setminus K},
\end{equation}
 whose weak Euler-Lagrange equation is given by
\begin{equation}\label{eq2} \begin{cases} u\in L^{1,p}(\Om\setminus K),&
u=g\quad\mbox{on}\quad\partial_D\Om\setminus K,\\ 
\salt
\displaystyle\int_{\Om\setminus K}f_\xi(x,\nabla u)\cdot\nabla \var\,\dx=0 & 
\forall\var\in L^{1,p}(\Om\setminus K),\quad\var =0\quad\mbox{on}\quad\partial_D\Om\setminus K.
\end{cases}
\end{equation}
By well-known existence results for nonlinear elliptic equations involving strictly monotone 
operators (see e.g. Lions \cite{lions}), one can easily see that (\ref{eq2}) 
has a unique solution in the sense that the gradient is always unique.

From now on, given $K\in\mathcal K_m(\overline \Om)$ and $u\in L^{1,p}(\Om\setminus K)$, we set
\begin{equation}\label{enemin}
E(u,K):=\int_{\Om\setminus K}f(x,\nabla u)\dx\,+\,\mathcal H^1(K).
\end{equation}

\begin{definition}\label{defunimin}
Let $g\in W^{1,p}(\Om)$ and let $m$ be a positive integer. We say that a pair $(u,K)$, with $K\in\mathcal K_m(\overline \Om)$, 
 $u\in L^{1,p}(\Om\setminus K)$ and $u=g$ on $\partial_D\Om\setminus K$ is an unilateral minimum of (\ref{enemin}) if 
\begin{equation}\label{euhemn}
E(u,K)\leq E(v,H)\end{equation}
 %$$\int_{\Om\setminus K}f(x,\nabla u)\dx\,+\,\mathcal H^1(K)\,
%\leq\,\int_{\Om\setminus H}f(x,\nabla v)\dx\,+\,\mathcal H^1(H),$$
 among all $H\in\mathcal K_m(\overline \Om)$, $H\supset K$ and 
 $v\in L^{1,p}(\Om\setminus H)$ with $v=g$ on $\partial_D\Om\setminus H$.
\end{definition}  

%\begin{remark}\label{equidefmin}
%{\rm Note that it is not restrictive to assume that (\ref{euhemn}) holds 
% only for every $J\in\mathcal K_m(\overline \Om)$, 
%$J\supset K$ with $(J\setminus K)\cap\partial_N\Om=\emptyset$. 
% In fact, for every $H\in\mathcal K_m(\overline \Om)$, 
%$H\supset K$, we consider the set $\tilde J = H\setminus ((H\setminus K \cap) \partial_N\Om)$.
%Let $C_i$ all the connected component of $\tilde J$ 
%such that $C_i\cap K \neq \emptyset$ for every $i$. 
%We set $J$ as the union of such
%connected components. Note that $J\in\mathcal K_m(\overline \Om)$, $K\subseteq J\subseteq H$, so we have 
%\begin{eqnarray*} \int_{\Om\setminus K}f(x,\nabla u)\dx\,+\,\mathcal H^1(K\cup \partial_N\Om)\,
%&\leq&\int_{\Om\setminus J}f(x,\nabla v)\dx\,+\,\mathcal H^1(J\cup \partial_N\Om)\\
%&\leq&\int_{\Om\setminus J}f(x,\nabla v)\dx\,+\,\mathcal H^1(H\cup \partial_N\Om),
%\end{eqnarray*}
%that is, (\ref{euhemn}) holds.
%}
%\end{remark}

\subsection{Hausdorff convergence}\label{haus}
 We recall here the {\it Hausdorff distance} 
between two closed sets $K_1$ and $K_2$ defined by 
$$
d_H(K_1,K_2):=\max\graffe{\displaystyle\sup_{x\in K_1}{\rm dist}\,(x,K_2)\,,\,
\displaystyle\sup_{x\in K_2}{\rm dist}\,(x,K_1)},$$
 with the conventions ${\rm dist}
\,(x,\emptyset)={\rm diam}\,(\Om)$ and $\sup\emptyset =0$, so that
$$
d_H(\emptyset\,,K)=
\begin{cases}
0 &\mbox{ if }K=\emptyset, \\
{\rm diam}\,(\Om)&\mbox{ if }K\neq\emptyset . 
\end{cases}
$$ 
Let $(K_h)$ be a sequence of compact subsets of $\overline\Om $. We say that $(K_h)$ 
converges to $K$ in the {\it Hausdorff metric} if $d_H(K_h\,,K)$ converges to $0$. 
It is well-known (see e.g., \cite[Blaschke's Selection Theorem]{falc}) that 
${\mathcal K}(\overline\Om)$ and ${\mathcal K}_m(\overline\Om)$ are
 compact with respect to the  Hausdorff convergence. Moreover, using Go\l\c{a}b theorem 
on the lower semicontinuity of the one-dimensional Hausdorff measure, we have 
 also that ${\mathcal K}^\lambda_m(\overline\Om)$ is compact with respect to the  Hausdorff convergence.
\vskip .1truecm
The following Lemma is proved  in \cite{dmro1}.
\begin{lemma}\label{dmtoa1}
Let $U$ be a bounded connected open subset of $\R^2$ with Lipschitz continuous boundary. 
Let $K$ be a closed connected subset of $\overline U$. Let $\lambda>0$ and let $(K_h)\subset {\mathcal K}^\lambda_m(\overline U)$ 
be a sequence which converges to $K$ in the Hausdorff metric. Then there exists a sequence $(H_h)$ of 
closed connected subsets of $\overline U$ which converges to $K$ in the Hausdorff metric, with $K_h\subset H_h$ 
for every $h$ and ${\mathcal H}^1(H_h\setminus K_h)\to 0$.
\end{lemma}
\begin{lemma}\label{dmtoa}
 Let $U$ be a bounded connected open subset of $\R^2$ with Lipschitz continuous boundary 
 and let $(K_h)\subset {\mathcal K}^f_m(\overline U)$ 
be a sequence which converges to a compact set $K$ in the Hausdorff metric. Let $\Gamma$ be a compact subset of $\overline U$ with 
 a finite number of connected components. Then there exists a sequence $(H_h)\subset {\mathcal K}^f_m(\overline U)$  
  which converges to $K$ in the Hausdorff metric, with $K_h\subset H_h$ 
for every $h$,  ${\mathcal H}^1(H_h\setminus K_h)\to 0$ and such that any connected component of $H_h\cup\Gamma$ converges to 
 a connected component of $K\cup\Gamma$ in the Hausdorff metric.
\end{lemma}
The proof of this lemma follows the lines of \cite[Lemma 3.6]{dmro1}. 
Precisely, we apply Lemma \ref{dmtoa1} to every connected component $C$ of $K\cup\Gamma$ 
and the union of those connected components of $K_h\cup\Gamma$ 
whose limits in the Hausdorff metric are contained in $C$.
\vskip .1truecm
The following Lemma  proved  in \cite{dmro1} will also be useful in the proof of our main results.

\begin{lemma}\label{approx}
Let $p$ and $m$ be two positive integers. Let $(K_h)$ be a sequence in $\mathcal K^f_p(\overline\Om)$ which converges 
in the Hausdorff metric to $K\in \mathcal K^f_p(\overline \Om)$, and let $H\in\mathcal K^f_m(\overline\Om)$ with $H\supset K$. 
Then there exists a sequence $(H_h)\subset\mathcal K^f_m(\overline\Om)$ such that $H_h\to H$ in the Hausdorff metric, $K_h\subset H_h$, 
and ${\mathcal H}^1(H_h\setminus K_h)\to\mathcal H^1(H\setminus K)$.
\end{lemma}

\vskip .1truecm
In order to study the continuity of the solution $u$ of (\ref{eq1}) with respect to 
the variations of the compact set $K$, we should be able to compare two solutions defined in 
two different domains. This is why, throughout this paper, given a function $u\in L^{1,p}(\Om\setminus K)$, 
 we extend $\nabla u$ in $\Om$ by setting $\nabla u=0$ in $\Om\cap K$. 
\subsection{Capacity}\label{cap}
 Let $1<r<\infty$ and let $B$ be a bounded open set in $\R^2$. 
For every subset $E$ of $B$, the $(1,r)$-capacity of $E$ in $B$, denoted by $C_r(E,B)$, 
is defined as the infimum of $ \int_B|\nabla u|^r\,dx$ over the set of all functions 
$u\in W^{1,r}_0(B)$ such that $u\geq 1$ a.e. in a neighborhood of $E$. 
If $r>2$, then $C_r(E,B)>0$ for every nonempty set $E$. On the 
contrary, if $r=2$ there are nonempty sets $E$ with $C_r(E,B)=0$ (for 
instance, $C_r(\{x\},B)=0$ for every $x\in B$).

We say that a property $\mathcal P(x)$ holds $C_r$-{\it quasi everywhere} (abbreviated $C_r$-{\it q.e.})
 in a set $E$ if it holds for all $x\in E$ except a subset $N$ of $E$ with $C_r(N,B)=0$.
 We recall that the expression {\it almost everywhere} (abbreviated {\it a.e.}) refers, as usual,
 to the Lebesgue measure. 

A function $u:E\to\overline\R$ is said to be {\it quasi-continuous} 
if for every $\e$ there exists $A_\e\subset E$, with  $C_r(A_\e,B)<\e$, such that the restriction 
of $u$ to $E\setminus A_\e$ is continuous. If $r>2$ every quasi-continuous function is continuous, 
while for $r=2$ there are quasi-continuous functions that are not 
continuous. It is well known that, 
for any open subset $U$ of $\R^2$, any function $u\in
 L^{1,r}(U)$ has a {\it quasi-continuous representative} 
$\overline u:U\cup\partial_LU\to\R$ which satisfies
$$
\lim_{\rho\to 0^+}\intcauchy_{B_\rho (x)\cap U}|u(y)-\overline u(x)|\,dy=0
\quad\mbox{for $C_r$-q.e. }x\in U\cup\partial_LU,
$$ 
where $\partial_LU$ denotes the Lipschitz part of the boundary $\partial U$ of $U$.
 We recall that if $u_h$ converges to $u$  strongly in $W^{1,r}(U)$, then 
a subsequence of $\overline u_h$ converges to $\overline u$ pointwise $C_r$-q.e. on $U\cup\partial_LU$.
 To simplify the notation we shall always identify throughout the paper each function $u\in L^{1,r}(U)$ with
 its quasi-continuous representative $\overline u$.
  
For these and other properties on quasi-continuous representatives
  the reader is referred to 
\cite{EG}, \cite{hekima}, \cite{maz1}, \cite{ziem}.
\vskip .1truecm
The following lemma is proved in \cite[Lemma 4.1]{dmro1} for $p=2$. 
The case $p\neq 2$ can be proved in the same way.
\begin{lemma}\label{lemmaaux}
Let $(K_h)$ be a sequence in $\mathcal K(\overline\Om)$ 
which converges to a compact set $K$ in the Hausdorff 
metric. Let $u_h\in L^{1,p}(\Om\setminus K_h)$ be a 
sequence such that $u_h=0$ $C_p$-q.e. on $\partial_D\Om\setminus K_h$ and
 $(\nabla u_h)$ is bounded in $L^p(\Om,\R^2)$. 
Then, there exists a function $u\in L^{1,p}(\Om\setminus K)$ with $u=0$ $C_p$-q.e. on 
$\partial_D\Om\setminus K$ such that, up to a subsequence,  
$\nabla u_h$ converges weakly to $\nabla u$ in $L^p(A,\R^2)$ 
for every $A\subset\subset\Om\setminus K$. 
If, in addition, $|K_h|$ converges to $|K|$, then 
 $\nabla u_h$ converges weakly to $\nabla u$ in $L^p(\Om,\R^2)$.
\end{lemma}
The following three lemmas will be crucial in the proof of our main result.

\begin{lemma}\label{moscostante} 
Let $(K_h) \subset {\mathcal K}_1(\overline\Om)$ converging to 
a compact set $K$ in the Hausdorff metric.
Let $(v_h)$ be a sequence in $W^{1,q}(\Om)$ converging  weakly in
$W^{1,q}(\Om)$ to a function $v$, with $v_h = 0$ \mbox{$C_q$-q.e.\,in $K_h$.} Then
$v =0$ $C_q$-q.e. in $K$.  
\end{lemma} 
\begin{proof}  We consider an open ball $B$ containing $\overline\Om$ and we extend both 
functions $v_h$ and $v$ to functions still denoted respectively by $v_h$ and $v$ such that 
the two extensions belong two $W^{1,q}_0(B)$ and $v_h \rightharpoonup v$ in $W^{1,q}(B)$.
Let $w_h$ and $w$ be the solutions of the problems
\begin{eqnarray}
 w_h\in W^{1,q}_0(B\setminus K_h),\quad\quad &\Delta_q w_h =
\Delta_q v\quad\mbox{in}\quad B\setminus K_h,\label{gdm1}
\\
w\in W^{1,q}_0(B\setminus K),\quad\quad &\Delta_q w =
\Delta_q v\quad\mbox{in}\quad B\setminus K.\nonumber
\end{eqnarray}
Using a result on the stability of Dirichlet problems by Bucur and Trebeschi
\cite{bt} (see also \v Sver\'ak \cite{sver}
for the case $q=2$), we obtain that $w_h$ converges to $w$ strongly 
in $W^{1,q}_0(B)$. Taking $v_h-w_h$ as test function in (\ref{gdm1}), 
which is possible since $v_h-w_h\in W^{1,q}_0(B\setminus K_h)$ (see, e.g., \cite[Theorem
4.5]{hekima}),
we obtain
\begin{equation}
\langle \Delta_q w_h , v_h-w_h\rangle =
\langle \Delta_q v , v_h-w_h\rangle,
\label{gdm2}
\end{equation}
where $\langle\cdot,\cdot\rangle$ is the duality pairing between 
$W^{-1,p}(B)$ and $W^{1,q}_0(B)$.
Passing to the limit in (\ref{gdm2}) we obtain
$$
\langle \Delta_q w , v-w\rangle =
\langle \Delta_q v , v-w\rangle,
$$
which implies $v=w$ by the strict monotonicity of $-\Delta_q$. Since, 
by definition, 
$w=0$ $C_q$-q.e.\ in $K$, we conclude that $v=0$ $C_q$-q.e.\ in $K$.
\end{proof}

\begin{lemma}\label{moscostante1}
Let $(K_h) \subset {\mathcal K}_1(\overline\Om)$ converging 
to a compact set $K$ in the Hausdorff metric.
Let $(v_h)$ be a sequence in $W^{1,q}(\Om)$, converging  weakly in
$W^{1,q}(\Om)$ to a function $v$. Assume that every function 
$v_h$ is constant $C_q$-q.e.\ in $K_h$. Then  $v$ is constant $C_q$-q.e.\ in $K$.
\end{lemma}
\begin{proof}
 This is trivial if $K$ contains only a single point. If $K$ 
has more than one point, there exists $r>0$ such that ${\rm
diam}(K_h)>2r$ for $h$ large enough. Let us prove that
the constant values $c_h$ taken by $v_h$ on $K_h$
are bounded uniformly with respect to~$h$. To this aim
let us consider a point $x_h\in K_h$. Since ${\rm
diam}(K_h)>2r$, we have $K_h\setminus B(x_h,r)\neq\emptyset$, and
by connectedness
\begin{equation}\label{polar}
K_h\cap\partial B(x_h,\rho)\neq\emptyset\qquad \hbox{ for every }
0<\rho<r.
\end{equation}
As $v_h=c_h$ $C_q$-q.e.\ on $K_h$, by using polar coordinates we
deduce from (\ref{polar}) the Poincar\'e inequality
$$
\int_{B(x_h,r)} |v_h-c_h|^q\dx\le M r^q \int_{B(x_h,r)} |\nabla
v_h|^q\dx,
$$
where the constant $M$ is independent of $h$ and $r$. Since 
the sequence $v_h$ is bounded in $W^{1,q}(\Om),$  it follows that $c_h$ is bounded, and so it
converges (up to a subsequence) to some constant $c$. So, the sequence
$v_h-c_h$ converges weakly to $v-c$ in 
$W^{1,q}(\Om)$, and by Lemma \ref{moscostante}
we get that $v = c$ $C_q$-q.e.\ on $K$. 
\end{proof}
\section{Conjugates and their properties}\label{harm}
Let $R$ be the rotation on $\R^2$ defined by $R(y_1,y_2):=(-y_2,y_1)$.  The following proposition
 on the global construction of conjugates will 
be crucial in the proof of Theorem \ref{unistabi}. 
%\subsection{Global conjugates}\label{global}
\begin{proposition}\label{conjglob}
Let $K\in {\mathcal K}(\overline\Om)$ and let $u$ be a solution of the problem (\ref{eq1}). 
Assume that $\Om$ is simply connected. Then there exists a function $v\in W^{1,q}(\Om)$ such that
$\nabla v=Rf_\xi(x,\nabla u)1_{K^c}$ a.e. in $\Om$.  Moreover, $v$ is constant $C_q$-q.e. on each 
connected component of $K\cup\partial_N\Om.$
\end{proposition}
\begin{proof} 
Let $u$ be a solution of (\ref{eq1}). We consider the vector field 
$\Phi\in L^q(\Om,\R^2)$ defined by 
$$\Phi:=f_\xi(x,\nabla u)1_{K^c}.$$ We have that 
${\rm div}(\Phi)=0$ in $\mathcal D'(\Om)$; hence
${\rm rot}(R\Phi)=0$ in $\mathcal D'(\Om)$. 
As $\Om$ is simply connected and has a Lipschitz boundary, 
there exists $v\in W^{1,q}(\Om)$ such that  
$\nabla v=R\Phi\mbox{ a.e. on }\Om$. 
 
Let us now prove
that $v$ is constant $C_q$-q.e. on each connected component of $K\cup\partial_N\Om$ we
proceed as follows. Let $C$ be a connected component of $K\cup\partial_N\Om$
with $C_{1,q}(C)>0$ and let $\e >0$. We set
 $$C_\e:=\{x\in
\overline{\Om}\,\mbox{: } {\rm dist}\,(x\,,\,C)<\e\}
\quad\mbox{and}\quad K_\e:=(K\cup\partial_N\Om)\cup \overline C_\e.$$
Let $u_\e$ be
the solution of the problem (\ref{eq1}) in $\Om\setminus K_\e$. From Lemma \ref{lemmaaux}
 applied to $u_\e-g$ and by the monotonicity of $K_\e$, we have that
$\nabla u_\e$ converges (up to a subsequence) to $\nabla u^*$ weakly in $L^p(\Om\setminus K,\R^2)$ 
for some $u^*\in L^{1,p}(\Om\setminus K)$ with $u^*=g$ on $\partial_D\Om\setminus K$. 

We claim that $\nabla u^*=\nabla u$ a.e. in $\Om$. 
 Indeed, by reformulating the problem (\ref{eq2}) as a variational inequality in $\Om\setminus K_\e$
 and using Minty's lemma, we get 
$$
\displaystyle\int_{\Om\setminus K_\e}f_\xi(x,\nabla z)\cdot (\nabla z -\nabla u_\e)
\dx\geq 0\qquad\forall z\in L^{1,p}(\Om\setminus K_\e),\, z
 = g\mbox{ on }\partial_D\Om\setminus K_\e.
$$
Now, let $ z\in L^{1,p}(\Om\setminus K)$ with $z =g$ on $\partial_D\Om\setminus K$. 
By the monotonicity of $K_\e$, we have that 
 $z\in L^{1,p}(\Om\setminus K_\e)$ 
and  $z = g\mbox{ on }\partial_D\Om\setminus K_\e$. So, 
$$
\displaystyle\int_{\Om\setminus K_\e}f_\xi(x,\nabla z)\cdot (\nabla z -\nabla u_\e)\dx\geq 0.
$$
 Using the convention that $\nabla u_\e=0$ in $\Om\cap K_\e$ we obtain
\begin{equation}\label{eqaux3}
\displaystyle\int_{\Om\setminus K}f_\xi(x,\nabla z)\cdot (\nabla z -\nabla u_\e)\dx
\geq -\int_{K_\e\setminus K}f_\xi(x,\nabla z)\cdot \nabla z\dx.
\end{equation}
 Now, letting $\e\to 0$ in (\ref{eqaux3}) we obtain
$$
\displaystyle\int_{\Om\setminus K}f_\xi(x,\nabla z)\cdot (\nabla z -\nabla u^*)\dx\geq 0 
\qquad\forall z\in L^{1,p}(\Om\setminus K),\, z
 = g\mbox{ on }\partial_D\Om\setminus K.
$$
which, using again Minty's lemma is equivalent to
$$
\displaystyle\int_{\Om\setminus K}f_\xi(x,\nabla u^*)\cdot\nabla\var\dx=0 
\qquad\forall\var\in L^{1,p}(\Om\setminus K),\, \var
 = 0\mbox{ on }\partial_D\Om\setminus K.
$$
By the uniqueness of solution of (\ref{eq2}) in $\Om\setminus K$, we get that $\nabla u^*=\nabla u$.
 So, we have proved
that all the sequence $(\nabla u_\e)$  converges  to $\nabla u$ weakly in
$L^p(\Om,\R^2)$. On the other hand, one can see that
\begin{equation}\label{limmono}
\lim_{\e\to 0}\int_\Om [f_\xi(x,\nabla u_\e) -f_\xi(x,\nabla u)]\cdot (\nabla u-\nabla u_\e)\dx=0.
\end{equation}
Hence arguing as in \cite[Lemma 2.4]{dmebpo} (recall that $f_\xi(x,\cdot)$ is  strictly monotone),
 it follows that $\nabla u_\e$  converges strongly to $\nabla u$ in
$L^p(\Om,\R^2)$.\\
Now, from the first part of the proof, we consider a function $v_\e\in W^{1,q}(\Om)$ such that 
$\nabla v_\e= Rf_\xi(x,\nabla u_\e)1_{K^c_\e}$ a.e. in $\Om$. We can assume 
that $\int_\Om v_\e\dx=\int_\Om v\dx=0$. So, by Poincar\'e inequality 
we obtain that $v_\e$ converges 
strongly to $v$ in $W^{1,q}(\Om)$. By construction $\nabla v_\e=0$ in $C_\e$ 
from which it follows that $v_\e$ is constant $C_q$-q.e. on $C_\e\cup\partial_LC_\e$. Hence 
$v_\e$ is constant $C_q$-q.e. on $C$. Since a subsequence of $v_\e$ 
converges to $v$ $C_q$-q.e. on $\overline\Om$, we conclude that $v$ 
is constant $C_q$-q.e. on $C$ and this completes the proof.
\end{proof}
\begin{definition}\label{defconj}
The function $v$ in Proposition \ref{conjglob} is called a conjugate of the function $u$.
\end{definition}
The following lemma is proved like in \cite[Theorem 4.3]{dmro1} for $f(x,\xi)=|\xi|^2$.
For the reader's convenience we will give here the proof 
of the present version.
\begin{lemma}\label{v-->usol}
Let $K\in\mathcal K_m(\overline\Om)$ and $u\in L^{1,p}(\Om\setminus K)$ 
with $u=g$ on $\partial_D\Om\setminus K$. 
Assume that there exists $v\in W^{1,q}(\Om)$ such that
$\nabla v=Rf_\xi(x,\nabla u)1_{K^c}$ a.e. in $\Om$ and that $v$ is 
constant $C_q$-q.e. on every connected component of $K\cup\partial_N\Om$.
 Then $u$ is solution of (\ref{eq2}).
\end{lemma}

\begin{proof}  Let $C^1,\ldots,C^l$ be the connected components of 
 $K\cup\partial_N\Om$. Since $v=c^i$ $C_q$-q.e on $C^i$, by 
\cite[Theorem 4.5]{hekima}, we can approximate $v$ strongly in $W^{1,q}(\Om)$ 
by a sequence of functions $v_n\in C_c^\infty(\R^2)$ that are constant in a 
suitable neighborhood $V_n^i$ of $C^i$. Let $\var\in L^{1,p}(\Om\setminus K)$ with
 $\var=0$ on $\partial_D\Om\setminus K$ and let $\var_n\in W^{1,p}_0(\Om\setminus K)$ 
such that $\var_n=\var$ in $\Om\setminus\bigcup_iV^i_n$. Then we have that 
\begin{equation}\label{eq31}
\int_{\Om}R\nabla v_n\nabla\var\dx=\int_{\Om\setminus K}R\nabla v_n\nabla\var_ndx=0,
\end{equation}
where the last equality follows from the fact that the vector field 
$R\nabla v_n$ is divergence free. Then passing to the limit in (\ref{eq31}) for $n\to\infty$, we get 
$$\int_{\Om\setminus K}f_\xi(x,\nabla u)\nabla\var\dx=\int_{\Om}R\nabla v\nabla\var\dx=0
\quad\forall\var\in L^{1,p}(\Om\setminus K)\mbox{ with }\var=0\mbox{ on }\partial_D\Om\setminus K.$$
So $u$ is a solution of (\ref{eq2}). 
\end{proof}
%\subsection{Local conjugates}\label{local}
The following Lemma on the local construction of conjugates will 
be used in the proof of Theorem \ref{genecase}.
\begin{lemma}\label{locconj}
Let $K\in\mathcal K(\overline\Om)$ and let $u$ be a solution of (\ref{eq1}) in 
$\Om\setminus K$. Let $U$ be an open rectangle such that $U\cap\Om$ is 
a non empty simply connected set. Then there exists a function 
$v\in W^{1,q}(U\cap\Om)$ such that 
$\nabla v=Rf_\xi(x,\nabla u)1_{K^c}$ a.e. in $U\cap\Om$. 
Moreover, $v$ is constant $C_q$-q.e. on each connected component of $\overline{U}\cap(K\cup\partial_N\Om)$.
\end{lemma}
\begin{proof}
We note that $u$ is solution of
the following problem
$${\displaystyle
\min\left\{\int_{(U\cap\Om)\setminus K} f(x,\nabla w)\dx
\mbox{: }w\in L^{1,p}(U\cap\Om)\setminus K)\,\,\mbox{and}\,\,w=u\,\,\mbox{on}
\,\,\partial(U\cap\Om)\setminus K\right\}}.
$$
% So we get ${\rm div}\,(f_\xi(x,\nabla u))1_{K^c})=0$ in $\mathcal D'(U\cap\Om)$ 
%and hence ${\rm rot}\,(Rf_\xi(x,\nabla u))1_{K^c})=0$ in $\mathcal D'(U\cap\Om)$. 
Since $U\cap\Om$ is simply connected, we can apply Proposition \ref{conjglob} 
with $\Om$ replaced by $U\cap\Om$. So,  there exists a function 
$v\in  W^{1,q}(U\cap\Om)$ such that $\nabla v=Rf_\xi(x,\nabla u)1_{K^c}$ a.e. in $U\cap\Om$,
 with $v$ constant $C_q$-q.e. on each connected component of $\overline{U}\cap(K\cup\partial_N\Om)$.
\end{proof}
\section{The stability results relative to problem ($\ref{euhemn}$)}
In this section we give the  stability results relative to problem ($\ref{euhemn}$). 
First of all, we prove in the following proposition, the stability of problem ($\ref{eq1}$) under the condition that
 any connected component of $K_h\cup\partial_N\Om$ converges 
to a connected component of $K\cup\partial_N\Om$ in the Hausdorff metric.

\begin{proposition}\label{lemdisjoint}
Let $\Om $ be a simply connected and bounded 
open subset of $\R ^2$ with Lipschitz continuous boundary. 
Assume  that $\partial_N\Om$ has $M$ connected components.
Let $\lambda >0$ and let $(K_h)\subset {\mathcal K}^\lambda_m(\overline\Om)$ 
be a sequence which converges to a compact set $K$ 
in the Hausdorff metric. Let $(g_h)$ be a sequence in $ W^{1,p}(\Om)$ 
which converges to $g$ strongly in  $W^{1,p}(\Om)$. 
Let $u_h\in L^{1,p}(\Om\setminus K_h)$ and $u\in L^{1,p}(\Om\setminus K)$ be the solutions 
 of the minimization problem (\ref{eq1}) with boundary data $g_h$ and $g$ respectively. 
%Let $K_h^1,\ldots,K^l_h$ be the connected components of $K_h\cup\partial_N\Om$ 
%and let $K^1,\ldots,K^l$ be the connected components of $K\cup\partial_N\Om.$
Assume that any connected component of $K_h\cup\partial_N\Om$ converges 
to a connected component of $K\cup\partial_N\Om$ in the Hausdorff metric.
Then $\nabla u_h$ converges strongly to $\nabla u$ in $L^p(\Om,\R^2)$. 
\end{proposition}

\begin{proof}
% Let us prove that $\nabla u_h$ converges strongly to $\nabla u$ in $L^p(\Om,\R^2)$. 
By the growth assumptions (\ref{f1}) on the function $f$, we have 
that $\nabla u_h$ and $f_\xi(x,\nabla u_h)$ are bounded 
respectively in $L^p(\Om,\R^2)$ and in $L^q(\Om,\R^2)$. 
So, applying Lemma \ref{lemmaaux} to $u_h-g_h$, we obtain that $\nabla u_h$ converges (up to a subsequence)  to 
 $\nabla u^*$ weakly in $L^p(\Om,\R^2)$ for some function 
$u^*\in L^{1,p}(\Om\setminus K)$ and $u^*=g$ on $\partial_D\Om\setminus K$. 

On the other hand, 
there exists  a vector field $\Psi\in L^q(\Om,\R^2)$
 such that  $f_\xi(x,\nabla u_h)\rightharpoonup \Psi$ weakly in $L^q(\Om,\R^2)$. Let us prove that 
$\Psi=f_\xi(x,\nabla u^*)$ a.e. in $\Om$. Since $|K_h|=|K|=0$ it is sufficient to prove 
that for every open ball $B\subset\subset\Om\setminus K$, $\Psi=f_\xi(x,\nabla u^*)$ a.e. in $B$.  
Note that by the Hausdorff complementary convergence we have $B\subset\subset\Om\setminus K_h$ for $h$
large enough. 

We may assume that the mean values of $u_h$ and $u^*$ on $B$ are zero. Thus
the Poincar\'e inequality and the Rellich theorem imply that $u_h\to u^*$ strongly in $L^p(B)$. Let $z\in W^{1,p}(B)$
and $\var\in C^\infty_c(B)$ with $\var\ge 0$. For $h$ large
enough we have $B\subset\subset\Om\setminus K_h$, thus by the monotonicity of $f_\xi(x,\cdot)$ we
have
\begin{equation}\label{gdm30}
\int_B (f_\xi(x,\nabla z)-f_\xi(x,\nabla u_h))\cdot(\nabla z-\nabla
     u_h)\var\dx\ge0.
\end{equation}
 We have also
$$
\int_B f_\xi(x,\nabla u_h)\cdot\nabla (( z-
u_h)\,\var)\dx =0,
$$
which, together with (\ref{gdm30}), gives
\begin{equation}\label{gdm31}
\int_B f_\xi(x,\nabla z)\cdot\nabla (( z-
u_h)\var)\dx
     -\int_B (f_\xi(x,\nabla z)-f_\xi(x,\nabla u_h))\cdot\nabla \var
     \,(z-u_h)\dx\ge 0.
\end{equation}
We can pass to the limit in each term of (\ref{gdm31}) and we get
\begin{equation}\label{gdm32}
\int_B f_\xi(x,\nabla z)\cdot\nabla (( z-
u^*)\var)\dx
-\int_B (f_\xi(x,\nabla z)-\Psi)\cdot\nabla \var
     \,(z-u^*)\dx\ge0.
\end{equation}
As ${\rm div}\,\Psi=0$ in ${\mathcal D}'(B)$, we have
\begin{equation}\label{gdm33}
\int_B \Psi\cdot\nabla (( z-
u^*)\var)\dx =0.
\end{equation}
{}{}From (\ref{gdm32}) and (\ref{gdm33}) we obtain
$$
\int_B (f_\xi(x,\nabla z)-\Psi)\cdot(\nabla z-\nabla
     u^*)\,\var\dx\ge0.
$$
As $\var$ is arbitrary, we get
$(f_\xi(x,\nabla z)-\Psi)\cdot(\nabla z-\nabla
     u^*)\ge0$
     a.e.\ in $B$.
In particular, taking $z(x):= u^*(x)\pm\e \eta\cdot x$, with $\eta\in\R^2$ and
$\e>0$, we obtain
$\pm (f_\xi(x,\nabla  u^*\pm\e
\eta)-\Psi)\cdot\eta\ge0$ a.e.\ in $B$.
As $\e$ tends to zero we get
$(f_\xi(x,\nabla  u^*)-\Psi)\cdot\eta=0$ a.e.\ in $B$,
which implies that $f_\xi(x,\nabla  u^*)=\Psi$ a.e.\ in $B$ by the
arbitrariness of $\eta$.\\
So we have proved that $f_\xi(x,\nabla u_h)\rightharpoonup f_\xi(x,\nabla  u^*)$ weakly in $L^q(\Om,\R^2)$.
Now let us prove that $u^*$ is a solution of (\ref{eq1}) in $\Om\setminus K$. 

Now we use the assumption that  $K^i_h$ converges to $K^i$ for every $i$.
By Proposition \ref{conjglob} there exists $v_h\in W^{1,q}(\Om)$  such that 
 $\nabla v_h=Rf_\xi(x,\nabla u_h)$ a.e. in $\Om$ with $v_h$  
constant $C_q$-q.e. on each connected component of $K_h\cup\partial_N\Om$.
 Since $f_\xi(x,\nabla u_h)$ converges to $f_\xi(x,\nabla  u^*)$ 
weakly in $L^q(\Om,\R^2)$, there exists a function $v\in W^{1,q}(\Om)$
 such that $v_h\rightharpoonup v$ weakly in $W^{1,q}(\Om)$ 
and $\nabla v=Rf_\xi(x,\nabla u^*)$ a.e. in $\Om$.
  Moreover, by Lemma \ref{moscostante1} we get that $v$ is 
constant $C_q$-q.e. on $K^i$ for every $i$. So from Lemma \ref{v-->usol} it follows that 
$u^*$ is a solution of (\ref{eq1}) in $\Om\setminus K$ and hence 
Thus, $\nabla u^*=\nabla u$ a.e. in $\Om$. Therefore, all the sequence $\nabla u_h$ 
converges to $\nabla u$ weakly in $L^p(\Om,\R^2)$. 

Now let us prove that $\nabla u_h$ 
converges to $\nabla u$ strongly in $L^p(\Om,\R^2)$. First of all, by lower semicontinuity we have that
\begin{equation}\label{eq10}
\int_\Om f(x,\nabla u)\dx\leq\liminf_{h\to\infty}\int_\Om f(x,\nabla u_h)\dx.\end{equation}
By the  convexity of $f(x,\cdot)$ we have also that
\begin{equation}\label{eq11}
\int_\Om f(x,\nabla u)\dx\geq\int_\Om f(x,\nabla u_h)\dx+
\int_\Om f_\xi(x,\nabla u_h)\cdot(\nabla u-\nabla u_h)\dx.\end{equation}

Since
$$
\int_\Om f_\xi(x,\nabla u_h)\cdot(\nabla u_h-\nabla g_h)\dx =0
\quad\mbox{and}\quad\int_\Om f_\xi(x,\nabla u)\cdot(\nabla u-\nabla g)\dx =0,
$$
 it follows that
\begin{eqnarray*}
\lim_{h\to\infty}\int_\Om f_\xi(x,\nabla u_h)\cdot(\nabla u-\nabla u_h)\dx 
& =&\lim_{h\to\infty}\int_\Om f_\xi(x,\nabla u_h)\cdot(\nabla u-\nabla g_h)\dx\\
 & =& \int_\Om f_\xi(x,\nabla u)\cdot(\nabla u-\nabla g)\dx=0.\end{eqnarray*}
%$$
%\int_\Om(f_\xi(x,\nabla u_h)-f_\xi(x,\nabla u))\cdot(\nabla u_h-\nabla u)\dx=
%-\int_\Om f_\xi(x,\nabla u)\cdot\nabla u_h \dx-\int_\Om f_\xi(x,\nabla u_h)\cdot\nabla u\dx.$$ 
%Hence passing to the limit for $h\to\infty$ we get
%$$\lim_{h\to\infty}\int_\Om(f_\xi(x,\nabla u_h)-f_\xi(x,\nabla u))\cdot(\nabla u_h-\nabla u)\dx
%=-2\int_\Om f_\xi(x,\nabla u)\cdot\nabla u\dx=0.$$

Hence passing to the limit in (\ref{eq11}) we get
$$\int_\Om f(x,\nabla u)\dx\geq\limsup_{h\to\infty}\int_\Om f(x,\nabla u_h)\dx,$$
which together with (\ref{eq10}) implies 
\begin{equation}\label{eq12}
\lim_{h\to\infty}\int_\Om f(x,\nabla u_h)\dx=\int_\Om f(x,\nabla u)\dx.
\end{equation}
Since $\nabla u_h\rightharpoonup\nabla u$ weakly in $L^p(\Om,\R^2)$, using 
the  strict convexity of $f(x,\cdot)$, it follows from (\ref{eq12}) that $\nabla u_h$ 
converges to $\nabla u$ strongly in $L^p(\Om,\R^2)$.
\end{proof}
\vskip .2truecm
We are now in a position to prove the main results of the paper.
\subsection{The case $\Om$ simply connected}

\begin{theorem}\label{unistabi}
Let $\Om $ be a simply connected and bounded open subset of $\R ^2$ with Lipschitz continuous boundary. 
Assume  that $\partial_N\Om$ has $M$ connected components.
Let $\lambda >0$ and $(K_h)\subset {\mathcal K}^\lambda_m(\overline\Om)$ 
be a sequence which converges to a compact set $K$ in the Hausdorff metric. Let $(g_h)$ be a sequence in $ W^{1,p}(\Om)$ 
which converges to $g$ strongly in  $W^{1,p}(\Om)$. Let $u_h$ be such that $(u_h,K_h)$ is an unilateral 
minimum relative to $g_h$ of the functional $E$ defined in (\ref{enemin}) 
 and let $u\in L^{1,p}(\Om\setminus K)$ be the solution of the minimization problem (\ref{eq1}). 
Then $\nabla u_h$ converges strongly to $\nabla u$ in $L^p(\Om,\R^2)$. Moreover, the pair $(u,K)$ is an unilateral 
minimum of the functional $E$ relative to $g$.
\end{theorem}

\begin{proof}[The proof of Theorem \ref{unistabi}.] \text{} {\it Step 1.}
Let us prove that $\nabla u_h$ converges strongly to $\nabla u$ in $L^p(\Om,\R^2)$. 
Let $K^1_h,\ldots,K^{n_h}_h$
 be the connected components of $K_h\cup\partial_N\Om$. As by assumption $n_h\leq m+M$, passing to
a subsequence we can assume that $n_h=n$ for every $h$ and that, for every $i\in\{1,\ldots,n\}$,
$K_h^i$ converges to some compact connected set $K^i$ in the Hausdorff metric.

If  $K^i\cap K^j=\emptyset$ for every $i\neq j$, then  $K^1,\ldots,K^n$ are exactly the
connected components of $K\cup\partial_N\Om$. So, by Proposition \ref{lemdisjoint}
it follows that $\nabla u_h$ converges strongly to $\nabla u$ in $L^p(\Om,\R^2)$.

Now we remove the assumption that 
$K^i\cap K^j=\emptyset$ for every $i\neq j$. 

%Let $C^1,\ldots,C^l$ ($l\leq m$) be the connected components of $K\cup\partial_N\Om$
% such that $C^i \not\subseteq \partial_N\Om$ for every   $1\leq i \leq l$. 
%Fix $i\in\{1,\ldots,l\}$. There exists $I_i\subset\{1,\ldots,n\}$ such that  $C^i=\bigcup_{j\in I_i}K^j$.
% and for every $j\in I_i$, there exists
%$\sigma_j^i\in I_i$, $\sigma_j^i\neq j$ such that $K^{\sigma_j^i}\cap K^j\neq\emptyset$.
%We set $$K_h^{I_i}:=\bigcup_{j\in I_i}K_h^j.$$
%We have that $K_h^{I_i}$ converges to $C^i$ in the Hausdorff metric. 
%By Lemma \ref{dmtoa}, there exists a sequence of closed connected sets $H_h^i$ of $\overline\Om$ 
%which converges to $C^i$ in the Hausdorff metric and such that $K_h^{I_i}\subset H_h^i$ and
%$\mathcal H^1(H_h^i\setminus K_h^{I_i})\to 0$.
 % Now, we set $$H_h:=\bigcup_{i=1}^l H_h^i.$$
%Note that $(H_h)\subset \mathcal K_{m}(\overline \Om)$. 
%Since $K_h \cup \partial_N\Om =\bigcup_{i=1}^l K_h^{I_i}\cup\partial_N\Om$, we have
%$$K_h\cup\partial_N\Om\subset H_h\cup\partial_N\Om
%\quad\mbox{ and }\quad\mathcal H^1((H_h\cup\partial_N\Om)\setminus (K_h\cup\partial_N\Om))\to 0.$$ 

 Applying  Lemma \ref{dmtoa} for $U=\Om$ and $\Gamma=\partial_N\Om$,
 we obtain a sequence $(H_h)\subset {\mathcal K}^f_m(\overline\Om)$  
  which converges to $K$ in the Hausdorff metric, with $K_h\subset H_h$ 
for every $h$,  ${\mathcal H}^1(H_h\setminus K_h)\to 0$ 
and such that any connected component of $H_h\cup\partial_N\Om$ converges to 
 a connected component of $K\cup\partial_N\Om$ in the Hausdorff metric.

We consider now  the following minimization problem
\begin{equation}\label{minaux1}
\displaystyle\min_{w}\graffe{\int_{\Om\setminus H_h}f(x,\nabla w)\dx\,\mbox{: }w\in 
L^{1,p}(\Om\setminus H_h)\,,\quad w=g_h\,\mbox{ on }\,\partial_D\Om\setminus H_h}.
\end{equation} 
Let $w_h\in L^{1,p}(\Om\setminus H_h)$ be the solution of (\ref{minaux1}). 
%Note that, by construction,
%the connected components of $H_h$ converge in the Hausdorff distance to those  
%the connected components of $K\cup\partial_N\Om$ that are not connected components of $\partial_N\Om$. 
%Then, the connected components of $H_h \cup \partial_N\Om$
%converge to the connected components of $K\cup\partial_N\Om$.
From Proposition \ref{lemdisjoint}, it
follows that $\nabla w_h$ converges  to $\nabla u$ strongly in $L^p(\Om,\R^2)$.
 Now using the fact that the pair $(u_h,K_h)$ is a unilateral minimum of the functional (\ref{enemin}), we get that 
\begin{eqnarray}\label{unimin1}
\qquad\limsup_{h\to\infty}\int_{\Om\setminus K_h} f(x,\nabla u_h)\dx & \leq &
\lim_{h\to\infty}\int_{\Om\setminus H_h}f(x,\nabla w_h)\dx \,+\,
\lim_{h\to\infty}{\mathcal H}^1(H_h\setminus K_h)\\
\salt 
\nonumber &  =&
\int_{\Om\setminus K}f(x,\nabla u)\dx.\end{eqnarray}
 Hence, recalling that $\nabla u_h$ 
converges to $\nabla u^*$ weakly in   $L^p(\Om,\R^2)$, we obtain 
\begin{eqnarray*}
& &\int_{\Om\setminus K} f(x,\nabla u^*)\dx\,\,\, \leq \,\,\,
\liminf_{h\to\infty}\int_{\Om\setminus K_h} f(x,\nabla u_h)\dx \,\,\, \leq \,\,\,\\
& \leq &
\limsup_{h\to\infty}\int_{\Om\setminus K_h} f(x,\nabla u_h)\dx\,\,\, \leq \,\,\,
 \int_{\Om\setminus K} f(x,\nabla u)\dx,
\end{eqnarray*}
which implies (by the uniqueness of solution of (\ref{eq1}) 
in $\Om\setminus K$) that $\nabla u^*=\nabla u$ a.e. in $\Om$. So, all the sequence
$\nabla u_h$ converges  to $\nabla u$ weakly in $L^p(\Om,\R^2)$ and 
$$\lim_{h\to\infty}\int_\Om f(x,\nabla u_h)\dx=\int_\Om f(x,\nabla u)\dx.$$
 Since $\xi\to f(x,\xi)$ 
is strictly convex, it follows that $\nabla u_h$ converges 
strongly to $\nabla u$ in $L^p(\Om,\R^2)$ and this achieves the proof of Step 1.
\vskip .2truecm
{\it Step 2.} Let us prove that the pair $(u,K)$ is an unilateral 
 minimum of the functional $E$ relative \mbox{to $g$.} Let $H\in {\mathcal K}_{m}(\overline\Om)$ with $K\subset H$ 
and let $w\in L^{1,p}(\Om\setminus H)$ with $w=g$ on $\partial_D\Om\setminus H$. 
By Lemma \ref{approx}, there exists a sequence $(H_h)\subset {\mathcal K}_{m}(\overline\Om)$
 such that $H_h\to H$ in the Hausdorff metric, $K_h\subset H_h$, 
and ${\mathcal H}^1(H_h\setminus K_h)\to\mathcal H^1(H\setminus K)$. From Lemma \ref{dmtoa}, we have a sequence 
$(\tilde H_h)\subset {\mathcal K}_{m}(\overline\Om)$ which converges to $H$ in the Hausdorff metric and such that
  $H_h\subset\tilde H_h$,  
$\mathcal H^1(\tilde H_h\setminus H_h)\to 0$ and, every connected component of $H_h\cup\partial_N\Om$ converges in the Hausdorff metric to 
 a connected component of $H\cup\partial_N\Om$. Let 
  $z_h\in L^{1,p}(\Om\setminus\tilde H_h)$ and $z\in L^{1,p}(\Om\setminus H)$ be the solutions of (\ref{eq1}) 
 with boundary data $g_h$ and $g$ respectively. From Proposition \ref{lemdisjoint} it follows that 
 $\nabla z_h\to \nabla z$ strongly in $L^p(\Om,\R^2)$.

Now using the fact that the pair $(u_h,K_h)$ is a unilateral minimum of the functional (\ref{enemin}), we get that
\begin{equation}\label{unimin6}
\int_{\Om\setminus K_h}f(x,\nabla u_h)\dx\,
\leq\,\int_{\Om\setminus \tilde H_h}f(x,\nabla z_h)\dx\,+
\,{\mathcal H}^1(\tilde H_h\setminus K_h).\end{equation}
 So, passing to the limit in (\ref{unimin6}) and using the 
fact that $\nabla u_h\to\nabla u$ strongly in $L^p(\Om,\R^2)$ and $\nabla z_h\to \nabla z$ strongly in $L^p(\Om,\R^2)$,  we obtain 
\begin{eqnarray*}
& { } &\int_{\Om\setminus K}f(x,\nabla u)\dx \,\,\, = \,\,\,
\lim_{h\to\infty}\int_{\Om\setminus K_h}f(x,\nabla u_h)\dx \,\,\, \leq \,\,\,\\
 &\leq & \lim_{h\to\infty}\int_{\Om\setminus\tilde H_h}f(x,\nabla z_h)\dx\,+\,
\lim_{h\to\infty}\,{\mathcal H}^1(\tilde H_h\setminus K_h)\,\,\, \leq  \,\,\,\\
&\leq &\int_{\Om\setminus H}f(x,\nabla z)\dx\,+\,{\mathcal H}^1(H\setminus K)
 \,\,\, \leq  \,\,\, \int_{\Om\setminus H}f(x,\nabla w)\dx\,+\,{\mathcal H}^1(H\setminus K),
\end{eqnarray*}
 which gives Step 2 and achieves the proof of the theorem.
\end{proof}
%\begin{remark}\label{remunimin}
%{\rm Actually it is sufficient to require that the unilateral  minimality  holds for every 
%$H\in \mathcal K_{m}^{\lambda + \epsilon}(\overline\Om)$ 
%for a fixed $\epsilon>0$. 
%
%}
%\end{remark}
\subsection{The general case}
Here we remove the assumption that $\Om$ is simply connected 
and we prove the stability theorem below using the local conjugates
 in Lemma \ref{locconj}.
\begin{theorem}\label{genecase}
Let $\Om $ be a bounded connected open subset of $\R ^2$ with Lipschitz continuous boundary. 
Assume  that $\partial_N\Om$ has $M$ connected components.
Let $\lambda >0$ and $(K_h)\subset {\mathcal K}^\lambda_m(\overline\Om)$ 
be a sequence which converges to a compact set $K$ in the Hausdorff metric. Let $(g_h)$ be a sequence in $ W^{1,p}(\Om)$ 
which converges to $g$ strongly in  $W^{1,p}(\Om)$. Let $u_h$ be such that $(u_h,K_h)$ is an unilateral 
minimum relative to $g_h$ of the functional $E$ defined in (\ref{enemin}) 
 and let $u\in L^{1,p}(\Om\setminus K)$ be the solution of the minimization problem (\ref{eq1}). 
Then $\nabla u_h$ converges strongly to $\nabla u$ in $L^p(\Om,\R^2)$. Moreover, the pair $(u,K)$ is an unilateral 
minimum of the functional $E$ relative to $g$.
\end{theorem}

\begin{proof}[Proof of Theorem \ref{genecase}.]
\text{}\\ 
First of all let us  prove that $\nabla u_h$ converges strongly to $\nabla u$ in $L^p(\Om,\R^2)$. 
By the growth assumptions (\ref{f1}) on $f$, we have that  
 $\nabla u_h$ is bounded in $L^p(\R^2,\R^2)$. By Lemma 
 \ref{lemmaaux} applied to $u_h-g_h$, we have that $\nabla u_h$ 
converges (up to a subsequence) to $\nabla u^*$ weakly 
in $L^p(\R^2,\R^2)$ for some $u^*\in L^{1,p}(\Om\setminus K)$ with 
 $u^*=g$ on $\partial_D\Om\setminus K$. We claim that $\nabla u^*=\nabla u$ a.e. in $\Om$. 

To this aim, we fix $r>0$ 
such that the minimum of the diameters of the connected 
components of $\Om^c$ is equal to $3r$. Using the fact that $\Om$ has a Lipschitz continuous boundary, we may find
    two families of open rectangles $(Q_i)_{i=1}^n$ and  $(U_i)_{i=1}^n$ such that, 
 for every $i\in\{1,\ldots,n\}$,\, $Q_i\subset\subset U_i$,\, 
$Q_i\cap\Om\neq\emptyset$ and $U_i\cap\Om$ is a Lipschitz domain and,  
%$\{x_1,\ldots,x_n\}\subset\overline\Om$ be such that
$$
\overline\Om\subset\bigcup_{i=1}^n Q_i\quad\mbox{ and }\quad \max_{1\leq i\leq n}{\rm diam}(U_i)=2r.$$
We set $$\eta:=\min_{1\leq i\leq n}d(Q_i,\partial U_i).$$
 For every $i\in\{1,\ldots,n\}$, the number of  connected components $C$ of $\overline U_i\cap K_h$ which intersect $Q_i$ 
is less or equal to $m +  \lambda/\eta$. Indeed, 
if $C$ intersects $\partial U_i$, 
then $\mathcal H^1(C)\geq\eta$ and hence,
their number is at most  $\lambda/\eta$. If $C\cap \partial U_i=\emptyset$, 
then $C$ is a connected component of $K_h$,  and 
their number is less or equal to $m$. Similarly the number of  connected components  of $\overline U_i\cap \partial_N\Om$ which intersect $Q_i$ 
is less or equal to $M  +\mathcal H^1(\partial_N\Om)/\eta$.
 Let $K^{i,1}_h,\ldots,K^{i,k_h}_h$ be all the connected components 
of $\overline U_i\cap K_h$ which intersect $Q_i$. 
Since $k_h\leq  m + \lambda/\eta$, 
passing to a subsequence, we can assume that
 $k_h=k$ for every $h$.
We set $$K^i_h: =\bigcup_{j=1}^k K^{i,j}_h.$$
 
Up to a subsequence, we have that $K_h^i$ converges in the Hausdorff metric 
to some compact set $K^i \in {\mathcal K}^\lambda_k (\overline{U_i\cap\Om})$.
Let $\Gamma^i$ be the union of those connected components of  $\overline U_i\cap \partial_N\Om$ which intersect $Q_i$.
 By Lemma \ref{dmtoa} applied to $U=U_i\cap\Om$ and $\Gamma=\Gamma_i$, we get a sequence $(H^i_h)\subset {\mathcal K}^f_k (\overline{U_i\cap\Om})$  
  which converges to $K^i$ in the Hausdorff metric, with $K^i_h\subset H^i_h$ 
for every $h$,  ${\mathcal H}^1(H^i_h\setminus K^i_h)\to 0$ and such that any connected component of 
$H^i_h\cup \Gamma^i$ converges to 
 a connected component of $K^i\cup\Gamma^i$ in the Hausdorff metric. We set  
$$H_h: =\bigcup_{i=1}^n H^i_h.$$
Note that 
$$(H_h) \in {\mathcal K}_m ( \overline\Om),\quad  H_h\supset K_h,\quad {\mathcal H}^1(H_h\setminus K_h)\to 0$$
and $H_h$ converges in the Hausdorff metric to the compact set $\tilde K:=\bigcup_{i=1}^nK^i.$
Moreover it is easy to see that $ \tilde K = K $.

We consider now the minimization problem
\begin{equation}\label{newmin}
\displaystyle\min_{w}\graffe{\int_{\Om\setminus H_h}f(x,\nabla w)\dx\,\mbox{: }w\in 
L^{1,p}(\Om\setminus H_h)\,,\quad w=g_h\,\mbox{ on }\,\partial_D\Om\setminus H_h}.
\end{equation}
Let $\tilde u_h\in L^{1,p}(\Om\setminus H_h)$ be the solution
of problem (\ref{newmin}).
Applying Lemma \ref{lemmaaux} to $\tilde u_h-g_h$, we get that $\nabla\tilde u_h$ converges (up to a subsequence)  to 
  $\nabla\tilde u^*$ weakly in $L^p(\Om,\R^2)$,
for some function $\tilde u^*$ in $L^{1,p}(\Om\setminus K)$ with 
$\tilde u^*=g$ on $\partial_D\Om\setminus K$. As in the proof of Proposition \ref{lemdisjoint}, we have also that
$f_\xi(x,\nabla\tilde u_h)$ converges to $f_\xi(x,\nabla\tilde u^*)$ weakly 
 in $L^q(\Om,\R^2)$.

Let us prove that $\nabla\tilde u^*=\nabla u$ a.e. in $\Om$. By a localization argument, it is sufficient to prove that 
 the function $\tilde u^*$ satisfies:
 \begin{equation}\label{locar1}
\begin{cases}
\displaystyle\int_{(Q_i\cap\Om)\setminus K}f_\xi(x,\nabla\tilde u^*)\cdot\nabla\var\dx=0, &\\
\salt
\forall\var\in C^\infty_c(Q_i)\mbox{ with }\var=0
\mbox{ on }(Q_i\cap\partial_D\Om)\setminus K. &
\end{cases}
\end{equation}
 Let $\tilde H_h^i:=H_h\cap\overline{U_i\cap\Om}$ and $\tilde \Gamma^i:=\partial_N\Om\cap\overline{U_i}$. 
 Since the diameter of $U_i$ is strictly less than the minimum of the diameters of the connected components of $\Om^c$, we have that 
  the open set $U_i\cap\Om$ is simply connected. So, by Lemma \ref{locconj}, there exists a function $v^i_h\in W^{1,q}(U_i\cap\Om)$ 
such that $\nabla v^i_h=Rf_\xi(x,\nabla\tilde u_h)$ a.e. in $U_i\cap\Om$ 
and $v^i_h$ is constant $C_q$-q.e. on the connected components of $\tilde H_h^i\cup\tilde \Gamma^i$. Since 
$H_h^i\cup\Gamma^i\subset\tilde H_h^i\cup\tilde \Gamma^i$, we have that any  connected component of $H_h^i\cup\Gamma^i$ is contained in a 
 connected component of $\tilde H_h^i\cup\tilde \Gamma^i$. So we have also that $v^i_h$ is constant $C_q$-q.e. on the connected components of
 $H_h^i\cup\Gamma^i$.
 From the fact that 
$f_\xi(x,\nabla\tilde u_h)$ converges $f_\xi(x,\nabla\tilde u^*)$ weakly 
 in $L^q(\Om,\R^2)$, it follows that  $v^i_h$ converges weakly to some function $v^i$  in $W^{1,q}(U_i\cap\Om)$ 
 such that  $\nabla v^i=Rf_\xi(x,\nabla\tilde u^*)$ a.e. 
in $U_i\cap\Om$.  
Since any connected component of $H^i_h\cup\Gamma^i$ converges to 
 a connected component of $K^i\cup\Gamma^i$ in the Hausdorff metric and 
$v^i_h$ is constant $C_q$-q.e. on the connected components of  $H_h^i\cup\Gamma^i$, we get from Lemma \ref{moscostante1}  that 
$v^i$ is constant $C_q$-q.e. on every connected component 
of $K^i\cup\Gamma^i$. Now applying Lemma \ref{v-->usol} with $\Om$ 
 replaced by $Q_i\cap\Om$, we get that $\tilde u^*$ satisfies (\ref{locar1}). 
Therefore $\nabla\tilde u^*=\nabla u$ a.e. in $\Om$.
So, $\nabla\tilde u_h$ converges weakly to $\nabla u$ in $L^p(\Om,\R^2)$ 
and $f_\xi(x,\nabla\tilde u_h)$ converges to $f_\xi(x,\nabla u)$ weakly 
 in $L^q(\Om,\R^2)$. Thus, arguing as in the proof of Proposition \ref{lemdisjoint}, we get  that 
$\nabla\tilde u_h$ converges to  $\nabla u$ strongly in $L^p(\Om,\R^2)$.
 
Now, from the minimality of  the pair $(u_h,K_h)$, we have that 
\begin{equation}\label{unimin3}
\int_{\Om\setminus K_h}f(x,\nabla u_h)\dx\,
\leq\,\int_{\Om\setminus H_h}f(x,\nabla\tilde u_h)\dx\,+
\,{\mathcal H}^1(H_h\setminus K_h).\end{equation}
So, passing to the limit in (\ref{unimin3}) and using the 
fact that $\nabla u_h\rightharpoonup\nabla u^*$ 
weakly in $L^p(\Om,\R^2)$,  we obtain 
\begin{eqnarray*}
& { } &\int_{\Om\setminus K}f(x,\nabla u^*)\dx  \leq 
\liminf_{h\to\infty}\int_{\Om\setminus K_h}f(x,\nabla u_h)\dx\,\leq\,\\
&\leq &\lim_{h\to\infty}\int_{\Om\setminus H_h}f(x,\nabla\tilde u_h)\dx\dx\,+\,
\lim_{h\to\infty}\,{\mathcal H}^1(H_h\setminus K_h)\,=\,\int_{\Om\setminus K}f(x,\nabla u)\dx,
\end{eqnarray*}
 which implies (by the uniqueness of solution of (\ref{eq1}) in $\Om\setminus K$) that $\nabla u^*=\nabla u$ a.e. in $\Om$ and 
\begin{equation}\label{eq20}
\lim_{h\to\infty}\int_{\Om}f(x,\nabla u_h)\dx=\int_{\Om}f(x,\nabla u)\dx.\end{equation} 
Since $\nabla u_h\rightharpoonup\nabla u$ weakly in $L^p(\Om,\R^2)$, using 
the  strict convexity of $f(x,\cdot)$, it follows from (\ref{eq20}) that $\nabla u_h$ 
converges to $\nabla u$ strongly in $L^p(\Om,\R^2)$.  
This achieves the proof of the first part of the theorem.

Now let us prove that the pair $(u,K)$ is an unilateral 
 minimum of the functional $E$ relative to $g$. 
Let $H\in {\mathcal K}_{m}(\overline\Om)$ with $K\subset H$ 
and let $w\in L^{1,p}(\Om\setminus H)$ with $w=g$ on $\partial_D\Om\setminus H$. 
It is not restrictive to assume that $H\in {\mathcal K}^f_{m}(\overline\Om)$. 
By Lemma \ref{approx}, there exists a sequence $(H_h)\subset {\mathcal K}_{m}(\overline\Om)$
 such that $H_h\to H$ in the Hausdorff metric, $K_h\subset H_h$, 
and $${\mathcal H}^1(H_h\setminus K_h)\to\mathcal H^1(H\setminus K).$$
 Since $(K_h)\subset\mathcal K^\lambda_m(\overline\Om)$ and $H\in {\mathcal K}^f_{m}(\overline\Om)$, we have that 
 $(H_h)\subset {\mathcal K}^{\lambda +\e}_{m}(\overline\Om)$ for some $\e>0$.
Arguing as in the first part of the proof, we can construct a sequence
$(\tilde H_h)\subset {\mathcal K}_{m}(\overline\Om)$ such that $H_h\subset\tilde H_h$,  
$\mathcal H^1(\tilde H_h\setminus H_h)\to 0$, 
and denoting  $z_h\in L^{1,p}(\Om\setminus\tilde H_h)$ and $z\in L^{1,p}(\Om\setminus H)$  the solutions of (\ref{eq1}) 
with boundary data $g_h$ and $g$ respectively, we get that 
$\nabla z_h\to \nabla z$ strongly in $L^p(\Om,\R^2)$.
Now we can achieve the proof as in Step 2 of Theorem \ref{unistabi}. 
\end{proof}

\vskip .5truecm
\centerline{\sc Acknowledgements}
 \vskip .1truecm
\noindent The authors wish to thank Gianni Dal Maso for having proposed 
the subject of this paper, and for many interesting discussions. This work is part of the
European Research Training Network ``Homogenization and Multiple
Scales'' under contract HPRN-2000-00109, and of the Research \
Project ``Calculus of Variations''
supported by SISSA and by
the Italian Ministry of Education, University, and Research.
%The research of F. Ebobisse is supported by a S.I.S.S.A postdoctoral fellowship.
 
\end{document}